\documentclass{amsart}
\usepackage{amsmath,amsthm,amssymb,amsfonts,color}

\usepackage{graphicx}
\usepackage{tikz}

\usepackage{epsfig}
\usepackage[pagebackref, colorlinks=true]{hyperref}

\hypersetup{urlcolor=blue, citecolor=blue}

\def\doi#1{   {\href{http://dx.doi.org/#1}
		{{\mdseries\ttfamily DOI}}}}

\usepackage{graphics}

  \newcommand{\ep}{\varepsilon}

\newcommand{\R}{\mathbb{R}}
\newcommand{\N}{\mathbb{N}}

\newcommand{\half}{\frac{1}{2}}

\newcommand{\pt}{\partial_t}

\newcommand{\beeq}{\begin{equation}}\newcommand{\eneq}{\end{equation}}

\newcommand{\supp}{\mathrm{supp}}

\newenvironment{prf}{\noindent {\bf Proof.} }{\endprf\par}
\def \endprf{\hfill  {\vrule height6pt width6pt depth0pt}\medskip}

\def\<{\langle}             \def\>{\rangle}
\def\({\left(}                 \def\){\right)}
\numberwithin{equation}{section}
\newtheorem{thm}{Theorem}[section]

\newtheorem{lem}[thm]{Lemma}

\newtheorem{defn}[thm]{Definition}
\newtheorem{rem}[thm]{Remark}
\newtheorem{asum}[thm]{Assumption}

\title
[Blow-up for NLW with spatial derivatives]
{Blow-up of solutions to semilinear wave equations with spatial derivatives
}

\author{Kerun Shao}
\address{School of Mathematical Sciences\\ Zhejiang University\\Hangzhou 310058, P. R. China}\email{shaokr@163.com}

\author{Hiroyuki Takamura}
\address{Mathematical Institute\\ Tohoku University\\ Aoba, Sendai 980-8578, Japan}\email{hiroyuki.takamura.a1@tohoku.ac.jp}

\author{Chengbo Wang$^*$}\thanks{* Corresponding author}
\address{School of Mathematical Sciences\\ Zhejiang University\\Hangzhou 310058, P. R. China}\email{wangcbo@zju.edu.cn}

\keywords{Glassey conjecture, Semilinear wave equation, Blow-up, Lifespan}

\subjclass[2010]{35L05, 35L71, 35B30,
	35B33, 
	35B44}

\date{\today}

\begin{document}
	\bibliographystyle{plain}
	\maketitle

	\begin{abstract}
For small-amplitude semilinear wave equations with power type nonlinearity on the first-order spatial derivative,
the expected sharp upper bound on the lifespan of solutions is obtained for both critical cases and subcritical cases, for all spatial dimensions $n>1$. It is achieved uniformly by constructing the integral equations, deriving the ordinary differential inequality system, and iteration argument. Combined with the former works, 
the sharp lifespan estimates for this problem are completely established, at least for the spherical symmetric case.
	\end{abstract}
	
	\section{Introduction}
	Let $n\ge 1$, $p>1$. Consider the Cauchy problem for  semilinear wave equations with sufficiently small initial data (of size $\ep >0$):
	\begin{equation}\label{eq-Glassey-p_x}
		\left\{
		\begin{aligned}
			&(\pt^2-\Delta)u(t,x)=|\nabla_x u(t,x)|^p\ &&, \ (t,x)\in(0,T)\times\R^n \ ,\\
			&u(0,x)=\ep f(x)\ ,\ \partial_t u(0,x)=\ep g(x)\ &&, \ x\in\R^n.
		\end{aligned}
		\right.
	\end{equation}
	

	Here, the initial functions $f, g$ are assumed to be sufficiently smooth with sufficient decay. The  lifespan, denoted by $T(\ep)$, is the maximal existence time of the solution to \eqref{eq-Glassey-p_x}. So, $T(\ep)=\infty$ means there exists a global solution.
	As an analog of the semilinear wave equation
	\begin{equation*}
		\left\{
		\begin{aligned}
			&(\pt^2-\Delta) u=|\partial_t u|^p\\
			&u(0)=\ep f\ ,\ \partial_t u(0)=\ep g\ ,
		\end{aligned}
		\right.
	\end{equation*} 
	it is conjectured that, in problem \eqref{eq-Glassey-p_x}, the critical power, denoted by $p_c(n)$, for the global existence v.s. blow-up is also given by $\frac{n+1}{n-1}$, for $n\geq2$, and $\infty$, for $n=1$; see Glassey \cite{MR0711440Review}.
	
	Heuristically,
the solutions 
could behave like free waves,
with energy of size $\ep$ and decay rate $(n-1)/2$,
for the time interval $[0,T)$, when


$$\int_0^T [(1+t)^{-(n-1)/2}\ep]^{p-1} dt\ll 1 \ .$$
This gives rise to the following expected sharp estimate for the  lifespan, 
$$\ln(T(\ep))\ep^{p-1}\sim 1$$
for critical powers $p=p_c(n)$, 
	and
$$(T(\ep))^{1-\frac{n-1}2(p-1)}\ep^{p-1}\sim 1$$
for subcritical powers $p\in(1,p_c(n))$.

	Before we state our results, let us review some background. Concerning the nonlinear term $|\partial_t u|^p$, the global existence part of the conjecture is verified for general initial data in dimension $n=2,3$ by Hidano-Tsutaya \cite{MR1386769} and Tzvetkov \cite{MR1637692} independently, with the previous radial result in Sideris \cite{MR0711440} for $n=3$. When dimension $n\geq4$, the global existence part is established for any $p\in (p_c(n),1+2/(n-2))$ with radial initial data in Hidano-Wang-Yokoyama \cite{MR2980460}. For sample choice of the initial data, the blow-up results for all spatial dimensions are verified; see Rammaha \cite{MR0879355} and the references therein for $n\geq2$, except critical cases for the even dimensions, and Zhou \cite{MR1845748} for $n\geq1$, together with the sharp upper bound on the lifespan for $p\in(1,p_c(n)]$, which is 
\beeq\label{eq-Glassey-p_x-bu}\varlimsup_{\ep\rightarrow 0^+}	\int_0^T [(1+t)^{-(n-1)/2}\ep]^{p-1} dt<\infty, \forall p\in (1, p_c(n)]\ .\eneq


 Deduced from the well-posed theory for $p\in(1,p_c(n)]$, the sharp lower bound on the lifespan is also obtained for all spatial dimensions, that is,
\beeq\label{eq-Glassey-p_x-wp}
\varliminf_{\ep\rightarrow 0^+}	\int_0^T [(1+t)^{-(n-1)/2}\ep]^{p-1} dt>0, \forall p\in (1, p_c(n)]\ ;\eneq
	see Hidano-Wang-Yokoyama \cite{MR2980460} and Fang-Wang \cite{MR3169752}, for $n\geq2$, and Kitamura-Morisawa-Takamura \cite{MR4558951}, for $n=1$. 
	
Turning to the spatial nonlinear term $|\nabla_x u|^p$, the well-posed theory developed for $|\partial_t u|^p$ 
applies also in this setting, which also gives 
the expected sharp lower bound estimate of the lifespan.

 Let us focus on the blow-up results and the upper bound estimate. For dimension one, the sharp lifespan estimates are verified recently in Sasaki-Takamatsu-Takamura \cite{MR4643157}. 
 In the following, we restrict ourselves to the high dimensional case $n\ge 2$.
 
 For $n=3$, the blow-up result is first established in Sideris \cite{MR0711440} for $p=p(3)=2$.
	For $n=2$, the blow-up result is verified by Schaeffer \cite{MR0829595} for $p=p_c(2)=3$. For $n\geq4$, Rammaha \cite{MR0879355} obtained the blow-up result for all of the critical and subcritical powers, except the critical cases for even dimensions. 
	Surprisingly enough, we observe that the 
	blow-up result for $p=p_c(n)$ with even $n\ge 4$ is left open.
	
	Turning to the upper bound estimates, we recall that
	 Sideris \cite{MR0711440} also obtained the upper bound estimate
	\begin{equation*}
		\varlimsup_{\ep\rightarrow 0^+}\ln(T(\ep))\ep^2<\infty \  ,  p=p(3)=2\ .
	\end{equation*}
However, we remark that this upper bound does not agree with the well-known lower bound, which is the classical John-Klainerman's almost global result, $T(\ep)\ge \exp(c\ep^{-1})$.
Later, Rammaha \cite{MR1338808} studied the upper bound for $n=2, 3$ and $p=2$.
	Therefore, the sharp upper bound on the lifespan is still largely left open, which is conjectured 
to agree with the lower bound 
\eqref{eq-Glassey-p_x-wp} in general, that is,	
for sample choice of the initial functions, we have
\begin{equation}\label{eq-sharpupperbound-crit}
		\varlimsup_{\ep\rightarrow 0^+}\ln(T(\ep))\ep^{p-1}<\infty \ ,
	\end{equation} 
for critical powers $p=p_c(n)$, and
	\begin{equation}\label{eq-sharpupperbound-sub}
		\varlimsup_{\ep\rightarrow 0^+}T(\ep)\ep^{\frac{2(p-1)}{2-(n-1)(p-1)}}<\infty \ ,
	\end{equation}
	for subcritical powers $p\in(1,p_c(n))$.

The aim of this paper is to fill the long standing gap concerning the
blow-up and upper  bound of the lifespan, for all $n\ge 2$ and $p\in (1, p_c(n)]$.
	
	\section{Preliminaries and main results}
	\begin{asum}\label{asum-initialdata}
		From now on, we make the following assumption:
		\begin{itemize}
			\item The initial data $f,g\in C_c^\infty(\R^n)$, the space of smooth functions on $\R^n$ with compact support, are radial and non-negative. 
			\item Let $R$ be the minimum number such that $$\supp\ f, \supp\ g \subset \{x\in\R^n : |x|\leq R\} \ .$$
			\item For the sake of brevity, we require that
			\begin{equation}\label{eq-asum-fasum}
				A_{f}:=\int_{[\frac{3R}{4},R]\times \R^{n-1}} f dx>0\ .
			\end{equation}
		\end{itemize}
	\end{asum}
	According to local existence theorems in \cite[Theorem 1.2, Theorem 1.3]{MR2980460} and \cite[Theorem 1.1]{MR3169752}, for sufficiently small $\ep$, we have a unique local (weak) solution to \eqref{eq-Glassey-p_x} belonging to $C_tH_{rad}^2(\R^n)\cap C_t^1H_{rad}^1(\R^n).$ Here, $H_{rad}^m(\R^n)$ consists of all the spherical  symmetric functions lying in the usual Sobolev space $H^m(\R^n)$. In view of the discussion above, we can give a definition of the lifespan.
	\begin{defn}[Lifespan]
		Assume $\ep$ is sufficiently small. Under the Assumption \ref{asum-initialdata}, the lifespan $T(\ep)$ is defined to be the supremum of $T>0$ such that equation \eqref{eq-Glassey-p_x} is well-posed on $[0,T)\times\R^n$ in the sense of theorems from \cite{MR2980460,MR3169752}.
	\end{defn}
	Now, we are ready to present the main results.
	\begin{thm}[Critical cases]\label{thm-main-critical}
		Let $n\geq2$ and $p=p_c(n)$. Under the Assumption \ref{asum-initialdata}, the lifespan $T(\ep)$ of the equation \eqref{eq-Glassey-p_x} satisfies 
		\begin{equation}\label{eq-mainthm-sharpupperbound-crit}
			\varlimsup_{\ep\rightarrow 0^+}\ln(T(\ep))\ep^{p-1}<\infty \ .
		\end{equation} 
	\end{thm}
	\begin{thm}[Subcritical cases]\label{thm-main-subcrit}
		Let $n\geq2$ and $p\in(1, p_c(n))$. Under the Assumption \ref{asum-initialdata}, the lifespan $T(\ep)$ of the equation \eqref{eq-Glassey-p_x} satisfies 
		\begin{equation}\label{eq-mainthm-sharpupperbound-sub}
			\varlimsup_{\ep\rightarrow 0^+}T(\ep)\ep^{\frac{2(p-1)}{2-(n-1)(p-1)}}<\infty \ .
		\end{equation}
	\end{thm}
	\begin{rem}
		Actually, precise constants of the right hand sides of \eqref{eq-mainthm-sharpupperbound-crit} and \eqref{eq-mainthm-sharpupperbound-sub} are given in the proof; see \eqref{eq-prf-mainthm-crit-perciseconst} and \eqref{eq-prf-mainthm-subcrit-perciseconst}.
	\end{rem}

	\section{Two lemmas for systems of ordinary differential inequalities}
	Before establishing the proof of main theorems, we give two lemmas for systems of ordinary differential inequalities, which play a key role in proving the theorems.
	
	\subsection{Lemma for the critical cases}
	\begin{lem}\label{lem-crit}
		Given constants $T_0>0$ and $p>1$, for all $A>0$, let $\mathcal{A}_{crit}(A)$ denote the set consisting of $H\in C^2[0,\mathcal{T}_H)$ with $H(0)=H'(0)=0$, and
		\begin{alignat}{2}
			H''(t)&\geq A(t+1)^{-1}\ &&,\  0\leq t< \mathcal{T}_H \ , \label{eq-lem-crit-ineql1}\\
			H''(t)&\geq (t+1)^{-(p+1)}\ln^{-(p-1)}(t+1)H^p(t)\ &&,\  T_0\leq t < \mathcal{T}_H  \label{eq-lem-crit-ineql2}\ .
		\end{alignat}
	Then, $T_{crit}(A):=\sup\{\mathcal{T}_H|H\in\mathcal{A}_{crit}(A)\}$ is finite. More precisely, we have \begin{equation}\label{eq-lem-crit-upperbound}
		\varlimsup_{A\rightarrow 0^+}\ln(T_{crit}(A))A^{p-1}<\infty \ .
	\end{equation}
	\end{lem}	
	\begin{rem}\label{rem-crit}
		By the following proof, the inequality \eqref{eq-lem-crit-upperbound} can be specified to
		\begin{equation*}
			\varlimsup_{A\rightarrow 0^+}\ln(T_{crit}(A))A^{p-1}\leq2^{p-1}\frac{\max\{p,2(p-1)\}}{(p-1)^{3}}p^\frac{2p-1}{p-1} \ .
		\end{equation*}
	\end{rem}
	
	 For dimension 3, this lemma appeared in Rammaha \cite[Lemma 2]{MR1338808}.
	
	\begin{prf}
		Both $H$ and $H'$ are positive on $(0,\mathcal{T}_H)$, and $H$ is a convex function, due to $H(0)=H'(0)=0$ and the inequality \eqref{eq-lem-crit-ineql1}.
		\subsubsection*{Step 1}(Improvements on the logarithm term)
		We temporarily assume $\mathcal{T}_H=\infty$, and establish the improvement on the logarithm term.
		Integrating \eqref{eq-lem-crit-ineql1} twice from 0 to $t$, we have 
		\begin{equation}\label{eq-prf-crit-time1}
			H(t)\geq \frac{A}{2}(t+1)\ln(t+1) \ \mbox{for}\ t\geq e^2-1 \ .
		\end{equation}
		Let $T_1=\exp(\max\{2,\ln(T_0+1),2\sqrt{\frac{2}{p}},\frac{\sqrt{2p}}{p-1}\})-1$. Substituting the above inequality into \eqref{eq-lem-crit-ineql2}, we obtain
		\begin{equation}\label{eq-prf-crit-2ndprime1}
			H''(t)\geq(\frac{A}{2})^p(t+1)^{-1}\ln(t+1)\ \mbox{for} \ t\geq T_1 \ .
		\end{equation}
		Integrating \eqref{eq-prf-crit-2ndprime1} from $T_1$ to $t$, we find
		\begin{equation*}
			\begin{aligned}
				H'(t)&\geq\half(\frac{A}{2})^p\ln^2(t+1)\left(1-\left(\frac{\ln(T_1+1)}{\ln(t+1)}\right)^2\right)  \\
				&\geq\half(\frac{A}{2})^p(1-p^{-1})\ln^2(t+1)\ \mbox{for} \ t\geq \exp(p\ln(T_1+1))-1=:\tilde{T_1} \ ,
			\end{aligned}
		\end{equation*}
		as $p>1$. Then, we have
		\begin{equation*}
			H(t)\geq \half(\frac{A}{2})^p(1-p^{-1})\int_{\tilde{T_1}}^{t}\ln^2(s+1)ds \ \mbox{for} \ t\geq\tilde{T_1}\ .
		\end{equation*}
		Since
		\begin{equation*}
			\int_{\tilde{T_1}}^{t}\ln^2(s+1)ds=(t+1)\ln^2(t+1)-(\tilde{T_1}+1)\ln^2(\tilde{T_1}+1)-2\int_{\tilde{T_1}}^{t}\ln(s+1)ds \ ,
		\end{equation*}
		it follows, for all $t\geq T_2:=\exp(\sqrt{2p}\ln(\tilde{T_1}+1))-1$,
		\begin{equation}\label{eq-prf-crit-T1}
			\ln(t+1)\geq2p\max\{2,\frac{p}{p-1}\}\geq4p\ ,
		\end{equation}
		and then
		\begin{equation*}
			\begin{gathered}
				2p(\tilde{T_1}+1)\ln^2(\tilde{T_1}+1)\leq(t+1)\ln^2(t+1)\ ,  \\
				4p\int_{\tilde{T_1}}^{t}\ln(s+1)ds\leq4p(t+1)\ln(t+1)\leq(t+1)\ln^2(t+1) \ .
			\end{gathered}
		\end{equation*}
		Consequently, we obtain 
		\begin{equation}\label{eq-prf-crit-time2}
			H(t)\geq \half(\frac{A}{2})^p(1-p^{-1})^2(t+1)\ln^2(t+1)\ \mbox{for}\ t\geq T_2\ ,  
		\end{equation}
		which suggests that we can use the iterate method to improve the inequality \eqref{eq-prf-crit-time1}.
		
		Now we claim that 
		\begin{equation}\label{eq-prf-crit-timek}
			H(t)\geq C_k(t+1)\ln^{q_k}(t+1)\ \mbox{for}\ t\geq T_k\ , k=1, 2, \cdots \ , 
		\end{equation}
		where the sequences $\{C_k\}, \{q_k\}$, and $\{T_k\}$ satisfy
		\begin{equation*}
			\left\{
				\begin{aligned}
					C_{k+1}&=\frac{(p-1)^3}{\max\{p,2(p-1)\}}\frac{C_k^p}{p^{k+1}}\ &&,\ C_1=\frac{A}{2}\ , \\
					q_{k+1}&=pq_k-p+2\ &&,\ q_1=1\ , \\
					T_{k+1}&=\exp((2p)^{\frac{1}{q_{k+1}}}p\ln(T_k+1))-1\ &&. 
				\end{aligned}
			\right.
		\end{equation*}
		Let us verify this claim. By simple calculation, we deduce that
		\begin{align}
			&q_{k+1}=(p^k+p-2)/(p-1)\leq p^k\max\{\frac{2}{p},\frac{1}{p-1}\} \ ,\nonumber\\
			&\ln(T_{k+1}+1)=(2p)^{\sum_{i=1}^{k}\frac{1}{q_{i+1}}}p^k\ln(T_1+1) \ \label{eq-prf-crit-T_k}.
		\end{align}
		Hence, substituting \eqref{eq-prf-crit-timek} into \eqref{eq-lem-crit-ineql2} and integrating the result from $T_k$ to $t$, one finds
		\begin{equation*}
			H'(t)\geq(1-p^{-1})\frac{(C_k)^p}{q_{k+1}}\ln^{q_{k+1}}(t+1)\ \mbox{for} \  t\geq \exp(p\ln(T_k+1))-1=:\tilde{T_k} \ .
		\end{equation*}
		Noticing that $\ln(T_{k+1}+1)=(2p)^\frac{1}{q_{k+1}}\ln(\tilde{T_k}+1)$, and
		\begin{equation*}
			\begin{aligned}
				&\int_{\tilde{T_k}}^{t}\ln^{q_{k+1}}(s+1)ds\\
				=&(t+1)\ln^{q_{k+1}}(t+1)-(\tilde{T_k}+1)\ln^{q_{k+1}}(\tilde{T_k}+1)-q_{k+1}\int_{\tilde{T_k}}^{t}\ln^{q_{k+1}-1}(s+1)ds \ ,
			\end{aligned}
		\end{equation*}
	 	then for all $t\geq T_{k+1}$, we have
		\begin{equation}\label{eq-prf-crit-T_kesti}
			\begin{aligned}
				\ln(t+1)&\geq(2p)^{\sum_{i=1}^{k}\frac{1}{q_{i+1}}}p^k\ln(T_1+1)\\
				&\geq(2p)^{\half+\sum_{i=1}^{k}\frac{1}{q_{i+1}}}p^k\max\{\frac{2}{p},\frac{1}{p-1}\}\ \\
				&\geq 2pq_{k+1} \ , 
			\end{aligned}
		\end{equation}
		and
		\begin{equation*}
			\begin{aligned}
				2p(\tilde{T_k}+1)\ln^{q_{k+1}}(\tilde{T_k}+1)&\leq(t+1)\ln^{q_{k+1}}(t+1)\ ,  \\
				2pq_{k+1}\int_{\tilde{T_k}}^{t}\ln^{q_{k+1}-1}(s+1)ds\leq2pq_{k+1}(t+1)&\ln^{q_{k+1}-1}(t+1)\leq(t+1)\ln^{q_{k+1}}(t+1)\ .
			\end{aligned}
		\end{equation*}
		Now, it follows that, for $t\geq T_{k+1}$, 
		\begin{equation*}
			\begin{aligned}
				H(t)&\geq (1-p^{-1})^2\frac{C_k^p}{q_{k+1}}(t+1)\ln^{q_{k+1}}(t+1) \\
				&=\frac{(p-1)^3}{p^2}\frac{C_k^p}{p^k+p-2}(t+1)\ln^{q_{k+1}}(t+1)\\
				&\geq\frac{(p-1)^3}{\max\{p,2(p-1)\}}\frac{C_k^p}{p^{k+1}}(t+1)\ln^{q_{k+1}}(t+1) \\
				&=C_{k+1}(t+1)\ln^{q_{k+1}}(t+1)\ ,
			\end{aligned}
		\end{equation*}
		which completes the proof of the claim \eqref{eq-prf-crit-timek}.
		
		By an elementary calculation, we have 
		\begin{equation*}
			C_{k+1}=B_{k+1}(\tilde{C}_{crit}A)^{p^k} \ ,
		\end{equation*} 
		where
		\begin{align*}
			B_{k+1}=p^{\frac{(k+2)(p-1)+1}{(p-1)^2}}\left[\frac{(p-1)^3}{\max\{p,2(p-1)\}}\right]^{-\frac{1}{p-1}}\ , \nonumber\\
			\tilde{C}_{crit}=\half \left[\frac{(p-1)^3}{\max\{p,2(p-1)\}}\right]^{\frac{1}{p-1}}p^{-\frac{2p-1}{(p-1)^2}}\ . 
		\end{align*}
		Therefore, 
		\begin{equation}\label{eq-prf-crit-finaltimek}
			\begin{aligned}
				H(t)&\geq B_{k+1}(\tilde{C}_{crit}A)^{p^k}(t+1)\ln^{q_{k+1}}(t+1)\\
				&=B_{k+1}(\tilde{C}_{crit}A\ln^{\frac{1}{p-1}}(t+1))^{p^k}(t+1)\ln^{\frac{p-2}{p-1}}(t+1) \  \mbox{for}  \ t\geq T_{k+1}\ .
			\end{aligned}
		\end{equation}
		
		\subsubsection*{Step 2}(Improvement from the logarithm term to the power term)
		Now we are going to show that $\varlimsup_{A\rightarrow 0^+}\tilde{C}_{crit}A\ln^{\frac{1}{p-1}}(T_{crit}(A)+1)\leq1$. Seeking a contradiction, suppose that there exists a constant $\delta\in(0,1)$ and a sequence $A_j\rightarrow0^+$ such that 
		\begin{equation}
			\tilde{C}_{crit}A_j\ln^{\frac{1}{p-1}}(T_{crit}(A_j)+1)>1+\delta \ ,
		\end{equation}
		which means that there exits a function sequence $H_j\in\mathcal{A}_{crit}(A_j)$ such that
		\begin{equation}\label{eq-prf-crit-contra}
			\tilde{C}_{crit}A_j\ln^{\frac{1}{p-1}}(\mathcal{T}_{H_j}+1)>1+\delta \ .
		\end{equation}
		
		Because the formula \eqref{eq-prf-crit-T_k} means that the growth of $T_{k}$ is mainly an exponential function. For each sufficiently small $A_j$, we can always select a greatest $k_j\in \N$ such that 
		\begin{equation}\label{eq-prf-crit-ATlink}
			1+\frac{\delta}{4}\geq\tilde{C}_{crit}A_j\ln^{\frac{1}{p-1}}(T_{k_j+1}+1)>1+\frac{\delta}{8} \ ,
		\end{equation}
		and notice that $\lim_{j\rightarrow\infty}k_j=\infty$. Thus, for $t\in [T_{k_j+1},\mathcal{T}_{H_j})$, interpolating \eqref{eq-prf-crit-finaltimek} with \eqref{eq-lem-crit-ineql2}, we have
		\begin{equation}\label{eq-prf-crit-mu(t)}
			\begin{aligned}
				H_j''(t)&\geq (t+1)^{-(p+1)}\ln^{-(p-1)}(t+1)H_j(t)H_j^{p-1}(t) \\
				&\geq B_{k_j+1}^{p-1}(\tilde{C}_{crit}A)^{(p-1)p^{k_j}}(t+1)^{-2}\ln^{p^{k_j}-1}(t+1)H_j(t) \\
				&=:\mu(t)H_j(t) \ .
			\end{aligned}
		\end{equation}
		In the following content, the subscripts of $k_j$ and $H_j$ will be omitted if there is no risk of misunderstanding.
		For all sufficiently small $A_j$, because of \eqref{eq-prf-crit-ATlink} and \eqref{eq-prf-crit-T_k}, it follows
		\begin{equation*}
			T_{k+1}^2\mu(T_{k+1})\geq\left(\frac{T_{k+1}}{T_{k+1}+1}\right)^2B_{k+1}^{p-1}(1+\frac{\delta}{8})^{p^k(p-1)}\ln^{-1}(T_{k+1}+1)>1 \ .
		\end{equation*}
		Multiplying \eqref{eq-prf-crit-mu(t)} by $H'(t)$ and integrating the resulting inequality from $T_{k+1}$ to $t$, one has
		\begin{equation*}
			\begin{aligned}
				(H'(t))^2&\geq (H'(T_{k+1}))^2+\int_{T_{k+1}}^{t}\mu(s)[H^2]'(s)ds \\
				&\geq (H'(T_{k+1}))^2+(T_{k+1}^2\mu(T_{k+1}))^{-1}\int_{T_{k+1}}^{t}\mu(s)[H^2]'(s)ds \\
				&\geq\frac{\mu(t)H^2(t)}{T_{k+1}^2\mu(T_{k+1})}+(H'(T_{k+1}))^2-T_{k+1}^{-2}H^2(T_{k+1}) \\
				&\geq\frac{\mu(t)H^2(t)}{T_{k+1}^2\mu(T_{k+1})} \\
				&\geq(t+1)^{-2}\(\frac{\ln(t+1)}{\ln(T_{k+1}+1)}\)^{p^k-1}H^2(t) \ ,
			\end{aligned}
		\end{equation*}
		where one uses $\mu'(t)<0$ for $t\geq T_{k+1}$, $\forall k\geq1$, by \eqref{eq-prf-crit-T_k}, in the third inequality, and the convexity of $H(t)$, hence $tH'(t)>H(t), \forall t>0$, in the fourth inequality. So we have 
		\begin{equation*}
			H'(t)\geq(t+1)^{-1}\(\frac{\ln(t+1)}{\ln(T_{k+1}+1)}\)^{\frac{p^k-1}{2}}H(t) \ \mbox{for} \ t\in [T_{k+1},\mathcal{T}_H) \ .
		\end{equation*}
		Let $\bar{T}=\exp((1+\delta/4)^{p-1}\ln(T_{k+1}+1))-1$. By \eqref{eq-prf-crit-ATlink} and \eqref{eq-prf-crit-contra}, we have $\bar{T}<\mathcal{T}_H$.
		Then, it follows that, for $t\in[\bar{T},\mathcal{T}_H)$, 
		\begin{equation*}
			\begin{aligned}
				\ln\(\frac{H(t)}{H(T_{k+1})}\)&\geq\frac{2[\ln^{\frac{p^k+1}{2}}(t+1)-\ln^{\frac{p^k+1}{2}}(T_{k+1}+1)]}{(p^k+1)\ln^{\frac{p^k-1}{2}}(T_{k+1}+1)} \\
				&\geq\frac{2[1-(1+\frac{\delta}{4})^{-(p-1)\frac{p^k+1}{2}}]}{p^k+1}(1+\frac{\delta}{4})^{(p-1)\frac{p^k-1}{2}}\ln(t+1) \ .
			\end{aligned}
		\end{equation*}
		Therefore, for all sufficiently small $A_j$ such that
		\begin{equation*}
			\frac{2[1-(1+\frac{\delta}{4})^{-(p-1)\frac{p^k+1}{2}}]}{p^k+1}(1+\frac{\delta}{4})^{(p-1)\frac{p^k-1}{2}}\geq\frac{4p}{p-1} \ ,
		\end{equation*}
		we obtain
		\begin{equation*}
			\ln\(\frac{H(t)}{H(T_{k+1})}\)\geq\frac{4p}{p-1}\ln(t+1) \ \mbox{for}\  t\in[\bar{T},\mathcal{T}_H) \ ,
		\end{equation*}
		that is, 
		\begin{equation}\label{eq-prf-crit-Htpowgrow}
			H(t)\geq H(T_{k+1})(t+1)^{\frac{4p}{p-1}} \ .
		\end{equation}
		\subsubsection*{Step 3}(Standard ordinary differential inequality argument)
		Once again, interpolating \eqref{eq-prf-crit-Htpowgrow} with \eqref{eq-lem-crit-ineql2}, we obtain
		\begin{equation}\label{eq-prf-crit-2ndprimefinal}
			\begin{aligned}
				H''(t)&\geq(t+1)^{-(p+1)}\ln^{-(p-1)}(t+1)H^{\frac{p-1}{2}}(t)H^{\frac{p+1}{2}}(t) \\
				&\geq H^{\frac{p-1}{2}}(T_{k+1})H^{\frac{p+1}{2}}(t) \ \mbox{for} \ t\in[\bar{T},\mathcal{T}_H) \ .
			\end{aligned}
		\end{equation}
		By the convexity of $H(t)$,
		\begin{equation*}
			\begin{aligned}
				H(t)&\geq H(\bar{T})+(t-\bar{T})H'(\bar{T})\\
				&\geq H(\bar{T})+\frac{t-\bar{T}}{\bar{T}}H(\bar{T})\\
				&=\frac{t}{\bar{T}}H(\bar{T}) \ \mbox{for} \ t\geq\bar{T} \ .
			\end{aligned}
		\end{equation*}
		So, multiplying \eqref{eq-prf-crit-2ndprimefinal} by $H'(t)$ and integrating from $\bar{T}$ to $t$, we find
		\begin{equation}\label{eq-prf-crit-1stprimefinal}
			\begin{aligned}
				\half (H'(t))^2&\geq\frac{2}{p+3}H^{\frac{p-1}{2}}(T_{k+1})[H^{\frac{p+3}{2}}(t)-H^{\frac{p+3}{2}}(\bar{T})]\\
				&\geq\frac{2}{p+3}H^{\frac{p-1}{2}}(T_{k+1})[1-2^{-\frac{p+3}{2}}]H^{\frac{p+3}{2}}(t)\ \mbox{for} \ t\in [2\bar{T},\mathcal{T}_H) \ ,
			\end{aligned}
		\end{equation}
		as long as $A_j$ is so small that
		\begin{equation*}
			\begin{aligned}
				\tilde{C}_{crit}A_j\ln^\frac{1}{p-1}(2\bar{T}+1)&\leq\tilde{C}_{crit}A_j[\ln(\bar{T}+1)+\ln2]^\frac{1}{p-1} \\
				&\leq\tilde{C}_{crit}A_j\left[(1+\frac{\delta}{4})^{p-1}\ln(\bar{T}+1)+\ln2\right]^\frac{1}{p-1} \\
				&<(1+\frac{\delta}{4})^2 \\
				&<1+\delta \\
				&<\tilde{C}_{crit}A_j\ln^\frac{1}{p-1}(\mathcal{T}_H+1) \ .
			\end{aligned}	
		\end{equation*}
		Consequently, dividing both sides of \eqref{eq-prf-crit-1stprimefinal} by $H^{\frac{p+3}{2}}(t)$, taking the square root, and integrating from $2\bar{T}$ to $t$, we deduce that
		\begin{equation*}
			\begin{aligned}
				\(\frac{4(1-2^{-\frac{p+3}{2}})}{p+3}\)^\half H(T_{k+1})^\frac{p-1}{4}(t-2\bar{T})&\leq\frac{4}{p-1}H^{-\frac{p-1}{4}}(2\bar{T})\\
				&\leq\frac{4}{p-1}H^{-\frac{p-1}{4}}(T_{k+1}) \ .
			\end{aligned}
		\end{equation*}
		Finally, when $A_j$ is sufficiently small, by \eqref{eq-prf-crit-finaltimek}, we obtain
		\begin{equation*}
			\mathcal{T}_H\leq2\bar{T}+\frac{4}{p-1}\(\frac{4(1-2^{-\frac{p+3}{2}})}{p+3}\)^{-\half}H^{-\frac{p-1}{2}}(T_{k+1})\leq4\bar{T} \ ,
		\end{equation*}
		and then,by \eqref{eq-prf-crit-ATlink},
		\begin{equation*}
			\begin{aligned}
				\tilde{C}_{crit}A_j\ln^{\frac{1}{p-1}}(\mathcal{T}_{H_j}+1)&\leq\tilde{C}_{crit}A_j(2\ln2+\ln(\bar{T}+1))^{\frac{1}{p-1}} \\
				&\leq(1+\frac{\delta}{4})\tilde{C}_{crit}A_j\ln^{\frac{1}{p-1}}(\bar{T}+1) \\
				&\leq(1+\frac{\delta}{4})^2 \\
				&<1+\delta \ ,
			\end{aligned}
		\end{equation*}
		which is contradictory to \eqref{eq-prf-crit-contra}.
	\end{prf}

	\subsection{Lemma for the subcritical cases}
	\begin{lem}\label{lem-subcrit}
		Given $n\geq1$, $T_0>0$, and $p\in(1,p_c(n))$, for all $A>0$, let $\mathcal{A}_{sub}(A)$ denote the set consisting of $H\in C^2[0,\mathcal{T}_H)$ with $H(0)=H'(0)=0$, and
		\begin{alignat}{2}
			H''(t)&\geq A &&,\  0\leq t< \mathcal{T}_H \ , \label{eq-lem-subcrit-ineql1}\\
			H''(t)&\geq (t+1)^{-\frac{n+3}{2}p+\frac{n+1}{2}}H^p(t)\ &&,\  T_0\leq t < \mathcal{T}_H  \label{eq-lem-subcrit-ineql2}\ .
		\end{alignat}
		Then, $T_{sub}(A):=\sup\{\mathcal{T}_H|H\in\mathcal{A}_{sub}(A)\}$ is finite.. More precisely, we have
		\begin{equation}\label{eq-lem-subcrit-upperbound}
			\varlimsup_{A\rightarrow 0^+}T_{sub}(A)A^{\frac{2(p-1)}{2-(n-1)(p-1)}}<\infty \ .
		\end{equation}
	\end{lem}

	\begin{rem}\label{rem-subcrit}
		By the following proof, the inequality \eqref{eq-lem-subcrit-upperbound} can be specified to
		\begin{equation*}
			\varlimsup_{A\rightarrow 0^+}T_{sub}(A)A^{\frac{2(p-1)}{2-(n-1)(p-1)}}\leq\(4^{3p-1}(b_0b_1)^2p^\frac{4p}{p-1}\)^{\frac{1}{2-(n-1)(p-1)}} \ .
		\end{equation*}
		Here, $b_0$ and $b_1$ are given in \eqref{eq-lem-subcrit-b0b1}.
	\end{rem}
	
	This lemma is verified in Rammaha \cite[Lemma 1]{MR1338808} for $n=p=2$, and in Haruyama-Takamura \cite{2404.06274} for dimension 1.
	
	
	\begin{prf}
	Both $H$ and $H'$ are positive on $(0,\mathcal{T}_H)$, due to the inequality \eqref{eq-lem-subcrit-ineql1}. Integrating \eqref{eq-lem-subcrit-ineql1} from 0 to $t$ twice, we have
	\begin{equation*}
		H(t)\geq\frac{A}{2}t^2\geq\frac{A}{8}(t+1)^2 \ \mbox{for} \ t\geq\max\{1, T_0\}=: T_1 \ .
	\end{equation*}
	
	At this time, we claim that 
	\begin{equation}\label{eq-prf-subcrit-timek}
		H(t)\geq C_k(t+1)^{q_k} \ \mbox{for} \ t\geq T_k \ ,\ k=1, 2, \cdots \ ,
	\end{equation}
	where the sequences $\{C_k\}, \{q_k\}$, and $\{T_k\}$ satisfy
	\begin{equation*}
		\left\{
		\begin{aligned}
			&q_{k+1}=p(q_k-\frac{n+3}{2})+\frac{n+5}{2} \ &&, \ q_1=2 \ , \\
			&T_{k+1}=2^{\frac{1}{q_{k+1}-1}+\frac{1}{q_{k+1}}}(T_k+1)-1 \ &&, \\
			&C_{k+1}=\frac{1}{4b_0b_1}\frac{C_k^p}{p^{2k}} \ &&, \ C_1=\frac{A}{8} \ .
		\end{aligned}
		\right.
	\end{equation*}
	Here, 
	\begin{equation}\label{eq-lem-subcrit-b0b1}
		\begin{gathered}
			b_0=\frac{1}{p-1}-\frac{n-1}{2}+\max\{0,\frac{n+3}{2}-\frac{1}{p-1}\}\cdot p^{-1} \ , \\
			b_1=\frac{1}{p-1}-\frac{n-1}{2}+\max\{0,\frac{n+1}{2}-\frac{1}{p-1}\}\cdot p^{-1} \ .
		\end{gathered}
	\end{equation}
	
	It is easy to verify this claim. We substitute \eqref{eq-prf-subcrit-timek} into \eqref{eq-lem-subcrit-ineql2} and integrate the result to have
	\begin{equation*}
		H'(t)\geq\half\frac{C_k^p}{q_{k+1}-1}(t+1)^{q_{k+1}-1} \ \mbox{for} \ t\geq2^{\frac{1}{q_{k+1}-1}}(T_k+1)-1 \ .
	\end{equation*}
	Thus, 
	\begin{equation*}
		H(t)\geq\frac{1}{4}\frac{C_k^p}{q_{k+1}(q_{k+1}-1)}(t+1)^{q_{k+1}} \ \mbox{for} \ t\geq2^{\frac{1}{q_{k+1}}+\frac{1}{q_{k+1}-1}}(T_k+1)-1 \ .
	\end{equation*}
	By a simple calculation, we have
	\begin{equation*}
		q_{k+1}=p^k(\frac{1}{p-1}-\frac{n-1}{2})+\frac{n+3}{2}-\frac{1}{p-1} \ . 
	\end{equation*}
	Because $p\in (1,p_c(n))$, it follows that
	\begin{equation*}
			\frac{q_{k+1}}{p^k}\leq b_0\ , \ \frac{q_{k+1}-1}{p^k}\leq b_1 \ ,
	\end{equation*}
	which completes the proof of the claim.
	
	By an elementary calculation, we have 
	\begin{equation*}
		C_{k+1}=D_{k+1}(\tilde{C}_{sub}A)^{p^k} \ ,
	\end{equation*}
	where
	\begin{align*}
		D_{k+1}&=(4b_0b_1)^\frac{1}{p-1}p^{\frac{2[k(p-1)+p]}{(p-1)^2}} \ , \\
		\tilde{C}_{sub}&=\frac{1}{8}(4b_0b_1)^{-\frac{1}{p-1}}p^{-\frac{2p}{(p-1)^2}} \ . 
	\end{align*}
	Therefore, it follows that
	\begin{equation*}
		H(t)\geq D_{k+1}(\tilde{C}_{sub}A(t+1)^{\frac{1}{p-1}-\frac{n-1}{2}})^{p^k}(t+1)^{\frac{n+3}{2}-\frac{1}{p-1}} \ \mbox{for} \ t\geq T_{k+1} \ .
	\end{equation*}
	Let 
	\begin{equation*}
		\tilde{T}=2^{\sum_{k=1}^{\infty}\frac{1}{q_{k+1}}+\frac{1}{q_{k+1}-1}}(T_1+1)-1<\infty \ ,
	\end{equation*}
	since $p\in (1,p_c(n))$ implies $\frac{1}{p-1}-\frac{n-1}{2}>0$.
	For each $A\in(0,\tilde{C}_{sub}^{-1}(\tilde{T}+1)^{\frac{n-1}{2}-\frac{1}{p-1}})$, if there were a $t_A\in(0,\mathcal{T}_H)$ such that 
	\begin{equation*}
		\tilde{C}_{sub}A(t_A+1)^{\frac{1}{p-1}-\frac{n-1}{2}}>1>\tilde{C}_{sub}A(\tilde{T}+1)^{\frac{1}{p-1}-\frac{n-1}{2}}\ ,
	\end{equation*}
	we would find
	\begin{equation}\label{eq-prf-subcrit-gotoinfty}
		H(t_A)\geq D_{k+1}(\tilde{C}_{sub}A(t_A+1)^{\frac{1}{p-1}-\frac{n-1}{2}})^{p^k}(t_A+1)^{\frac{n+3}{2}-\frac{1}{p-1}} \ , \ \forall k\in \N \ ,
	\end{equation} 
	which contradicts that $H(t_A)$ is finite. Thus, we conclude that
	\begin{equation*}
		\varlimsup_{A\rightarrow 0^+}\tilde{C}_{sub}A(T_{sub}(A)+1)^{\frac{1}{p-1}-\frac{n-1}{2}}\leq1\ .
	\end{equation*}
	
	\end{prf}

\section{Proof of the theorem \ref{thm-main-critical} and theorem \ref{thm-main-subcrit}}
	We will construct an integral equation $U(t)$ near the wave front and derive the system of ordinary differential inequalities about $U(t)$. Then, we use the lemmas verified above to show the upper bound on the lifespan. The selection of $U(t)$ has appeared in the former works, e.g., \cite{MR1338808,MR1845748}.

	\begin{prf}
	Because of the local wellposedness, we can say that there exists a unique weak solution $u\in C([0,T(\ep));H_{rad}^2(\R^n))\cap C^1([0,T(\ep));H_{rad}^1(\R^n))$ to the problem \eqref{eq-Glassey-p_x}. By finite speed of propagation, we deduce that for all $t>0$, $\supp \ u(t,\cdot) \subset \{x\in\R^n : |x|\leq t+R\}$; see details in \cite[Lemma 2.11]{2211.01594}. For $r\in\R$, let a linear operator $*:w\mapsto w^*$ be
	\begin{equation*}
		w^*(t,r)=\int_{\R^{n-1}}w(t,r,\tilde{x})d\tilde{x} \ .
	\end{equation*}
	The operator $*$ is defined for all admissible functions.
	Then, $u^*$ is a weak solution to the following equation:
	\begin{equation*}
		\left\{
		\begin{aligned}
			&\partial_t^2u^*-\partial_r^2u^*=(|\nabla_x u|^p)^* \ ,\\
			&u^*(0,r)=\ep f^*(r) \ , \ \partial_t u^*(0,r)=\ep g^*(r) \ .
		\end{aligned}
		\right.
	\end{equation*}
	Notice that $u^*\in C([0,T(\ep));H^2(\R))\cap C^1([0,T(\ep));H^1(\R)).$
	Since $n\geq2$, $p\in(1,p_c(n)],$ and $\nabla_x u\in C([0,T(\ep));H^1(\R^n))\cap C^1([0,T(\ep));L^2(\R^n))$, 
	by H\"older inequalities, and Sobolev embedding for $H^1(\R^2)$, we have for all $t_1, t_2\in[0,T(\ep))$
	\begin{equation*}
		\begin{aligned}
			&\left|\int_{\R}(|\nabla_x u|^p)^*(t_1,r)-(|\nabla_x u|^p)^*(t_2,r)dr\right| \\
			=&\left|\int_{\R^n}|\nabla_x u(t_1,x)|^p-|\nabla_x u(t_2,x)|^pdx\right| \\
			\leq&p\int_{\R^n}\max\{|\nabla_x u(t_1,x)|^{p-1},|\nabla_x u(t_2,x)|^{p-1}\}|\nabla_x u(t_1,x)-\nabla_x u(t_2,x)|dx \\
			\leq&\left\{\begin{aligned}
					&p\|\nabla_x u\|_{L_t^\infty L^2(\R^n)}^{p-1}\|\nabla_x u(t_1,\cdot)-\nabla_x u(t_2,\cdot)\|_{L^2(\R^n)}[\alpha_n(t_1+t_2+R)^n]^\frac{2-p}{2} \ &&, \ p\leq2 \ ;   \\ 
					&p\|\nabla_x u\|_{L_t^\infty H^1(\R^2)}^{p-1}\|\nabla_x u(t_1,\cdot)-\nabla_x u(t_2,\cdot)\|_{L^2(\R^2)} \ &&, \  p\in(2,p_c(2)] \ ,  
				\end{aligned}
				\right. \\
		\end{aligned}
	\end{equation*}
	where $\alpha_n$ is the volume of the n-dimensional unit ball. Thus, we obtain that $(|\nabla_x u|^p)^*\in C([0,T(\ep));L^1(\R))$. Let $u_\delta^*=u^*(t,\cdot)\ast\eta_\delta$, where $\{\eta_\delta(r)\}_{\delta>0}$ are the standard modifiers.
	Using d'Alembert's formula and taking the limit of $u_\delta^*$, we have
	\begin{equation}\label{eq-prf-mainthm-1dimwaveformula}
		u^*(t,r)=\ep\bar{u}(t,r)+\half\int_{0}^{t}\int_{r-t+\tau}^{r+t-\tau}(|\nabla_x u|^p)^*(\tau,\lambda)d\lambda d\tau  \ ,
	\end{equation}
	where
	\begin{equation*}
		\bar{u}(t,r)=\frac{f^*(r-t)+f^*(r+t)}{2}+\half\int_{r-t}^{r+t}g^*(s)ds \ .
	\end{equation*}
	
	Define 
	\begin{equation*}
		U(t)=\int_{0}^{t}(t-\tau)\int_{\tau+R_0}^{\tau+R}r^{-\beta}u^*(\tau,r)drd\tau \ .
	\end{equation*}
	Here $\beta\geq0$ and $R_0\in[0,R)$ are constants to be fixed later.
	Obviously, $U\in C^2[0,T(\ep))$ with $U(0)=U'(0)=0$, and 
	\begin{equation}\label{eq-prf-mainthm-U''}
		U''(t)=\int_{t+R_0}^{t+R}r^{-\beta}u^*(t,r)dr \ .
	\end{equation}
	Due to Assumption \ref{asum-initialdata} and \eqref{eq-prf-mainthm-1dimwaveformula}, we have $\bar{u}\geq0$, hence $u\geq0$. By H\"older's inequality and finite propagation speed of $u$,
	\begin{equation*}
		\begin{aligned}
			(|\nabla_x u|^p)^*(\tau,\lambda)&\geq\int_{\R^{n-1}}|\partial_{x_1}u|^p(\tau,\lambda,\tilde{x})d\tilde{x} \\
			&\geq[\alpha_{n-1}((\tau+R)^2-\lambda^2)^{\frac{n-1}{2}}]^{-(p-1)}\left|\int_{\R^{n-1}}\partial_{x_1}u(\tau,\lambda,\tilde{x})d\tilde{x}\right|^p\\
			&=[\alpha_{n-1}((\tau+R)^2-\lambda^2)^{\frac{n-1}{2}}]^{-(p-1)}|\partial_ru^*|^p(\tau,\lambda)\\
			&=:m(\tau,\lambda)|\partial_ru^*|^p(\tau,\lambda) \ .
		\end{aligned}
	\end{equation*}
	Substituting \eqref{eq-prf-mainthm-1dimwaveformula} into \eqref{eq-prf-mainthm-U''}, we obtain 
	\begin{equation}\label{eq-prf-mianthm-2ndprimeU}
		\begin{aligned}
			U''(t)&=\ep\bar{U}(t)+ \half\int_{t+R_0}^{t+R}r^{-\beta}\int_{0}^{t}\int_{r-t+\tau}^{r+t-\tau}(|\nabla_x u|^p)^*(\tau,\lambda)d\lambda d\tau dr \\
			&\geq\ep\bar{U}(t)+ \half\int_{t+R_0}^{t+R}r^{-\beta}\int_{0}^{t}\int_{r-t+\tau}^{r+t-\tau}m(\tau,\lambda)|\partial_ru^*|^p(\tau,\lambda)d\lambda d\tau dr \ ,
		\end{aligned}
	\end{equation}
	where
	\begin{equation*}
		\bar{U}(t)=\int_{t+R_0}^{t+R}r^{-\beta}\bar{u}(t,r)dr\ .
	\end{equation*}
	Thus, we have, for $t\geq0$,
	\begin{equation}\label{eq-prf-mainthm-U''iter1st}
		U''(t)\geq\ep(t+R)^{-\beta}\int_{R_0}^{R}\frac{f^*(r)}{2}dr=\ep(t+R)^{-\beta}A_{f} \ .
	\end{equation}
	Using the Lemma 3 in \cite{MR1338808}, one finds, for $t\geq R_1:=(R-R_0)/2$, 
	\begin{equation}\label{eq-prf-mainthm-nonlinearpartintest}
		\begin{aligned}
			&\int_{t+R_0}^{t+R}r^{-\beta}\int_{0}^{t}\int_{r-t+\tau}^{r+t-\tau}m(\tau,\lambda)|\partial_ru^*|^p(\tau,\lambda)d\lambda d\tau dr \\
			\geq(t+&R)^{-\beta-1}\int_{0}^{t}\int_{\tau+R_0}^{\tau+R}(t-\tau)(\lambda-\tau-R_0)m(\tau,\lambda)|\partial_ru^*|^p(\tau,\lambda)d\lambda d\tau . 
		\end{aligned}
	\end{equation}
	Although the Lemma 3 in \cite{MR1338808} requires $\beta>0$, it is still available and can be easily verified, when $\beta=0$.
	It follows from H\"older's inequality that
	\begin{equation}\label{eq-prf-mainthm-intp_ru^*J(t)}
		\begin{aligned}
			&\int_{0}^{t}\int_{\tau+R_0}^{\tau+R}(t-\tau)(\lambda-\tau-R_0)m(\tau,\lambda)|\partial_ru^*|^p(\tau,\lambda)d\lambda d\tau \\
			\geq&\left|\int_{0}^{t}\int_{\tau+R_0}^{\tau+R}(t-\tau)(\lambda-\tau-R_0)\lambda^{-\beta}\partial_ru^*(\tau,\lambda)d\lambda d\tau \right|^pJ^{-(p-1)}(t) \ , 
		\end{aligned}
	\end{equation}
	where
	\begin{equation*}
		\begin{aligned}
			J(t)&=\int_{0}^{t}\int_{\tau+R_0}^{\tau+R}(t-\tau)(\lambda-\tau-R_0)\lambda^{-\beta p'}m^{-\frac{1}{p-1}}(\tau,\lambda)d\lambda d\tau \\
			&=\alpha_{n-1}\int_{0}^{t}\int_{\tau+R_0}^{\tau+R}(t-\tau)(\lambda-\tau-R_0)\lambda^{-\beta p'}((\tau+R)^2-\lambda^2)^{\frac{n-1}{2}}d\lambda d\tau\ ,
		\end{aligned}
	\end{equation*}
	and $1/p+1/p'=1.$
	Through integration by part, the integral in the right hand side of \eqref{eq-prf-mainthm-intp_ru^*J(t)} becomes
	\begin{equation*}
		-\int_{0}^{t}\int_{\tau+R_0}^{\tau+R}(t-\tau)\lambda^{-\beta}(1-\beta\frac{\lambda-\tau-R_0}{\lambda})u^*(\tau,\lambda)d\lambda d\tau .
	\end{equation*}
	Hence, we have, as long as $1-2\beta R_1/R_0>0$,
	\begin{equation}\label{eq-prf-mianthm-backtoU(t)}
		\begin{aligned}
			&\left|\int_{0}^{t}\int_{\tau+R_0}^{\tau+R}(t-\tau)(\lambda-\tau-R_0)m(\tau,\lambda)|\partial_ru^*|^p(\tau,\lambda)d\lambda d\tau\right|^p \\
			\geq&\(1-2\beta\frac{R_1}{R_0}\)^pU^p(t) \ . 
		\end{aligned}
	\end{equation}
	As for the term $J(t)$,
	\begin{equation*}
		\begin{aligned}
			J(t)&\leq\alpha_{n-1}(R-R_0)^{\frac{n+3}{2}}2^{\frac{n-1}{2}}\int_{0}^{t}(t-\tau)(\tau+R)^{\frac{n-1}{2}}(\tau+R_0)^{-\beta p'}d\tau \\
			&\leq\alpha_{n-1}R_1^{\frac{n+3}{2}}2^{n+1}\(\frac{R_0}{R}\)^{-\beta p'}(t+R)\int_{0}^{t}(\tau+R)^{\frac{n-1}{2}-\beta p'}d\tau \ . 
		\end{aligned}
	\end{equation*}
	By a simple calculation, we find
	\begin{equation*}
			\int_{0}^{t}(\tau+R)^{\frac{n-1}{2}-\beta p'}d\tau\leq\left\{
			\begin{aligned}
				&\frac{R^{\frac{n-1}{2}-\beta p'+1}}{\beta p'-\frac{n-1}{2}-1} \ &&, \ \frac{n-1}{2}-\beta p'<-1 \ ; \\
				&\ln\frac{t+R}{R} \ &&, \ \frac{n-1}{2}-\beta p'=-1 \ ; \\
				&\frac{(t+R)^{\frac{n-1}{2}-\beta p'+1}}{\frac{n-1}{2}-\beta p'+1}\ &&, \ \frac{n-1}{2}-\beta p'>-1 \ .
			\end{aligned}
			\right.
	\end{equation*}
	Thus, it follows that
	\begin{equation}\label{eq-prf-mainthm-J(t)J_p(t)}
		J(t)\leq c_{n,\beta,p}\bar{J}_p(t) \ ,
	\end{equation}
	where
	\begin{equation*}
		c_{n,\beta,p}=\alpha_{n-1}R_1^{\frac{n+3}{2}}2^{n+1}\(\frac{R_0}{R}\)^{-\beta p'}\cdot\left\{
		\begin{aligned}
			&\frac{R^{\frac{n-1}{2}-\beta p'+1}}{\beta p'-\frac{n-1}{2}-1} \ &&, \ \frac{n-1}{2}-\beta p'<-1 \ ; \\
			&1 \ &&, \ \frac{n-1}{2}-\beta p'=-1 \ ; \\
			&\frac{1}{\frac{n-1}{2}-\beta p'+1}\ &&, \ \frac{n-1}{2}-\beta p'>-1 \ ,
		\end{aligned}
		\right.
	\end{equation*}
	and
	\begin{equation*}
		\bar{J}_p(t)=\left\{
		\begin{aligned}
			&t+R \ &&, \ \frac{n-1}{2}-\beta p'<-1 \ ; \\
			&(t+R)\ln\frac{t+R}{R} \ &&, \ \frac{n-1}{2}-\beta p'=-1 \ ; \\
			&(t+R)^{\frac{n-1}{2}-\beta p'+2}\ &&, \ \frac{n-1}{2}-\beta p'>-1 \ .
		\end{aligned}
		\right.
	\end{equation*}
	Combining \eqref{eq-prf-mianthm-2ndprimeU}, \eqref{eq-prf-mainthm-nonlinearpartintest}, \eqref{eq-prf-mainthm-intp_ru^*J(t)}, \eqref{eq-prf-mianthm-backtoU(t)}, and \eqref{eq-prf-mainthm-J(t)J_p(t)}, we obtain, for $t\geq R_1$,
	\begin{equation}\label{eq-prf-mainthm-U''iter2nd}
		U''(t)\geq\half c_{n,\beta,p}^{-(p-1)}\(1-2\beta\frac{R_1}{R_0}\)^p(t+R)^{-\beta-1}\bar{J}_p^{-(p-1)}(t)U^p(t)\ .
	\end{equation}
	
	Hence, if we want to improve the estimate for $U''$ by \eqref{eq-prf-mainthm-U''iter1st} and \eqref{eq-prf-mainthm-U''iter2nd}, we can integrate \eqref{eq-prf-mainthm-U''iter1st} twice, substitute the result into \eqref{eq-prf-mainthm-U''iter2nd}, and hope the power of the term $t+R$ increases, that is,
	\begin{equation*}
		\left\{
		\begin{aligned}
			&-(p-1)+(2-\beta)p-\beta-1>-\beta \ &&, \ 1<p\leq\frac{n+1}{n+1-2\beta} \ ;\\
			&-(\frac{n-1}{2}-\beta p'+2)(p-1)+(2-\beta)p-\beta-1>-\beta \ &&, \ p>\frac{n+1}{n+1-2\beta}\ .
		\end{aligned}
		\right.
	\end{equation*}
	Equivalently,
	\begin{equation*}
		\left\{
		\begin{aligned}
			&(1-\beta)p>0 \ &&, \ 1<p\leq\frac{n+1}{n+1-2\beta} \ ;\\
			&p<p_c(n) \ &&, \ p>\frac{n+1}{n+1-2\beta}\ .
		\end{aligned}
		\right.
	\end{equation*}
	To sum up, for the subcritical cases, one expects to deduce improvements made in the term $t+R$, as long as $1-2\beta R_1/R_0>0$, $A_f>0$, and $\beta \in [0,1)$. As for the critical cases, improvements can only be made in the logarithm term by taking $\beta=1$.
	
	\subsection{Critical cases}
		Choose $\beta=1$ and $R_0=3R/4$, hence $1-2\beta R_1/R_0=2/3$. By \eqref{eq-asum-fasum}, we have $A_f>0.$
		According to \eqref{eq-prf-mainthm-U''iter1st} and \eqref{eq-prf-mainthm-U''iter2nd}, we have
		\begin{equation*}
			\left\{
			\begin{aligned}
				&U''(t)\geq\ep A_{f}(t+R)^{-1} \ &&, \ t\geq0 \ , \\
				&U''(t)\geq (2c_{n,1,p_c}^{-1})^{(p_c-1)}3^{-p_c}(t+R)^{-(p_c+1)}\ln^{-(p_c-1)}(\frac{t+R}{R})U^{p_c}(t) \ &&, \ t\geq R_1=\frac{R}{8} \ .
			\end{aligned}
			\right.
		\end{equation*}
	Set $\tilde{U}(t)=2\cdot3^{-\frac{p_c}{p_c-1}}c_{n,1,p_c}^{-1}R^{-1}U(Rt)$. It follows that $\tilde{U}(t)$ satisfies the inequalities that
	\begin{equation*}
		\left\{
		\begin{aligned}
			&\tilde{U}''(t)\geq2\cdot3^{-\frac{p_c}{p_c-1}} \ep A_{f}c_{n,1,p_c}^{-1}(t+1)^{-1}\ &&,\  t\geq0 \ ,\\
			&\tilde{U}''(t)\geq (t+1)^{-(p_c+1)}\ln^{-(p_c-1)}(t+1)\tilde{U}^{p_c}(t)\ &&,\  t\geq\frac{R_1}{R}=\frac{1}{8}  \ .
		\end{aligned}
		\right.
	\end{equation*}
	Therefore, by Remark \ref{rem-crit}, we deduce the upper bound estimate of the lifespan, which is 
	\begin{equation}\label{eq-prf-mainthm-crit-perciseconst}
		\begin{aligned}
			&\varlimsup_{\ep\rightarrow 0^+}\ln(\frac{T^*(\ep)}{R})\ep^{p_c-1} \\
			\leq&4^{-\frac{n+1}{n-1}}R^\frac{n+3}{n-1}\max\{n+1,4\}(n+1)^\frac{n+3}{2}(n-1)^{-\frac{n-1}{2}}(A_f^{-1}\alpha_{n-1})^\frac{2}{n+1} \ .
		\end{aligned}
	\end{equation}

	\subsection{Subcritical cases}
	We take $\beta=0$ to handle all the subcritical powers. Choose $R_0=3R/4$. By \eqref{eq-asum-fasum}, we have $A_f>0$. Thus, according to \eqref{eq-prf-mainthm-U''iter1st} and \eqref{eq-prf-mainthm-U''iter2nd}, we obtain
	\begin{equation*}
		\left\{
		\begin{aligned}
			&U''(t)\geq\ep A_{f} \ &&, \ t\geq0 \ , \\
			&U''(t)\geq\half c_{n,0,p}^{-(p-1)}(t+R)^{-\frac{n+3}{2}p+\frac{n+1}{2}}U^p(t) \ &&, \ t\geq \frac{R}{8} \ .
		\end{aligned}
		\right.
	\end{equation*}
	Set $\tilde{U}(t)=2^{-\frac{1}{p-1}}c_{n,0,p}^{-1}R^{\frac{1}{p-1}-\frac{n+3}{2}}U(Rt)$. It follows that $\tilde{U}(t)$ satisfies the inequalities that
	\begin{equation*}
		\left\{
		\begin{aligned}
			&\tilde{U}''(t)\geq \ep A_{f}2^{-\frac{1}{p-1}}c_{n,0,p}^{-1}R^{\frac{1}{p-1}-\frac{n-1}{2}}\ &&,\  t\geq0 \ ,\\
			&\tilde{U}''(t)\geq (t+1)^{-\frac{n+3}{2}p+\frac{n+1}{2}}\tilde{U}^p(t)\ &&,\  t\geq\frac{1}{8}  \ .
		\end{aligned}
		\right.
	\end{equation*}
	Hence, by Remark \ref{rem-subcrit}, we deduce the upper bound estimate,
	\begin{equation}\label{eq-prf-mainthm-subcrit-perciseconst}
		\begin{aligned}
			&\varlimsup_{\ep\rightarrow 0^+}T^*(\ep)\ep^{\frac{2(p-1)}{2-(n-1)(p-1)}} \\
			\leq&\left[2^{-(n-1)p+n+5}(b_0b_1)^2p^\frac{4p}{p-1}\(R^\frac{n+3}{2}A_{f}^{-1}\frac{\alpha_{n-1}}{n+1}\)^{2(p-1)}\right]^\frac{1}{2-(n-1)(p-1)}\ ,
		\end{aligned}	
	\end{equation}
	where $b_0$ and $b_1$ are given in \eqref{eq-lem-subcrit-b0b1}.
		
	\end{prf}
	
\section*{Acknowledgement}

The first and the third authors were supported  by  NSFC 12141102.

The second author is partially supported by the Grant-in-Aid for Scientific
Research (A) (No. 22H00097), Japan Society for the Promotion of Science.

	
	\bibliographystyle{plain}
        \bibliography{blowup_reference}


		\end{document}